\documentstyle[amssymb,12pt]{amsart}

\newtheorem{lem}{Lemma}
\newtheorem{prop}{Proposition}

\theoremstyle{definition}

\theoremstyle{definition}
\newtheorem{thm}{Theorem}

\theoremstyle{remark}
\newtheorem{rem}{Remark}

\numberwithin{equation}{section}

\begin{document}

\newcommand{\thmref}[1]{Theorem~\ref{#1}}
\newcommand{\secref}[1]{Sect.~\ref{#1}}
\newcommand{\lemref}[1]{Lemma~\ref{#1}}
\newcommand{\propref}[1]{Proposition~\ref{#1}}
\newcommand{\corref}[1]{Corollary~\ref{#1}}
\newcommand{\conjref}[1]{Conjecture~\ref{#1}}
\newcommand{\remref}[1]{Remark~\ref{#1}}

\newcommand{\nc}{\newcommand}
\nc{\on}{\operatorname}
\nc{\W}{{\cal W}}
\nc{\pa}{\partial}
\nc{\sw}{{\frak s}{\frak l}}
\nc{\sun}{\widehat{\sw}_n}
\nc{\g}{{\frak g}}
\nc{\n}{{\frak n}}
\nc{\G}{\widehat{\g}}
\nc{\wh}{\widehat}
\nc{\wt}{\widetilde}
\nc{\GG}{{\bold G}}
\nc{\dd}{\wt{\delta}}
\nc{\NN}{{\bold N}}
\nc{\gL}{^L\!\g}
\nc{\k}{h^\vee}
\nc{\beq}{\begin{equation}}
\nc{\C}{{\Bbb C}}
\nc{\arr}{\rightarrow}
\nc{\lsln}{L\sw_n}
\nc{\F}{{\cal F}_n}
\nc{\Fq}{{\cal F}_{n,q}}
\nc{\Z}{{\Bbb Z}}
\nc{\M}{{\cal M}_n}
\nc{\Mq}{{\cal M}_{n,q}}
\nc{\La}{\Lambda}
\nc{\al}{\alpha}
\nc{\rr}{\wt{r}}
\nc{\de}{\wt{\delta}}
\nc{\ot}{\otimes}
\nc{\nb}{\nabla}
\nc{\vf}{\varphi}
\nc{\la}{\lambda}
\nc{\pol}{\frac{1}{2}}
\nc{\RN}{{\cal R}_{n,q}}
\nc{\ol}{\overline}
\nc{\ZN}{\Z/N\Z}
\nc{\MM}{{\bold M}}
\nc{\ep}{\epsilon}
\nc{\D}{{\frak d}}
\nc{\vff}{\phi \left( \frac{t}{s} \right)}
\nc{\vpp}{\varphi \left( \frac{w}{z} \right)}
\nc{\muu}{{\bold m}}

\title[Drinfeld-Sokolov reduction for difference
operators I.]{Drinfeld-Sokolov reduction for difference operators and
deformations of ${\cal W}$--algebras I. \\ The case of Virasoro
algebra}

\author[E. Frenkel]{E. Frenkel$^*$}\thanks{$^*$ Packard Fellow}

\address{Department of Mathematics, Harvard University, Cambridge,
MA 02138, USA}

\author{N. Reshetikhin}

\address{Department of Mathematics, University of California, Berkeley, CA
94720, USA}

\author{M.A. Semenov-Tian-Shansky}

\address{Universit\'{e} de Bourgogne, Dijon, France and Steklov
Mathematical Institute, St. Petersburg, Russia}

\maketitle

\begin{abstract}
We propose a $q$--difference version of the Drin\-feld-Sokolov
reduction scheme, which gives us $q$--deformations of the classical
${\cal W}$--algebras by reduction from Poisson-Lie loop groups. We
consider in detail the case of $SL_2$. The nontrivial consistency
conditions fix the choice of the classical $r$-matrix defining the
Poisson-Lie structure on the loop group $LSL_2$, and this leads to a
new elliptic classical $r$--matrix. The reduced Poisson algebra
coincides with the deformation of the classical Virasoro algebra
previously defined in \cite{FR}. We also consider a discrete analogue
of this Poisson algebra. In the second part \cite{SS} the construction
is generalized to the case of an arbitrary semisimple Lie algebra.
\end{abstract}

\section{Introduction}

It is well-known that the space of ordinary differential operators of
the form $\pa^n + u_1 \pa^{n-2} + \ldots + u_{n-1}$ has a remarkable
Poisson structure, often called the (second) Adler-Gelfand-Dickey
bracket \cite{A,GD}. Drinfeld-Sokolov reduction \cite{DS} gives a
natural realization of this Poisson structure via the hamiltonian
reduction of the dual space to the affine Kac-Moody algebra
$\sun$. Drinfeld and Sokolov \cite{DS} have applied an analogous
reduction procedure to the dual space of the affinization $\G$ of an
arbitrary semisimple Lie algebra $\g$. The Poisson algebra ${\cal
W}(\g)$ of functionals on the corresponding reduced space is called
the classical ${\cal W}$--algebra. Thus, one can associate a classical
${\cal W}$--algebra to an arbitrary semisimple Lie algebra $\g$. In
particular, the classical ${\cal W}$--algebra associated to $\sw_2$ is
nothing but the classical Virasoro algebra, i.e. the Poisson algebra
of functionals on the dual space to the Virasoro algebra (see, e.g.,
\cite{FR}).

It is interesting that ${\cal W}(\g)$ admits another description as
the center of the universal enveloping algebra of an affine
algebra. More precisely, let $Z(\widehat{{\frak g}})_{-h^{\vee }}$ be
the center of a completion of the universal enveloping algebra
$U(\G)_{-\k}$ at the critical level $k=-h^{\vee }$ (minus the dual
Coxeter number). This center has a canonical Poisson structure. It was
conjectured by V.~Drinfeld and proved by B.~Feigin and E.~Frenkel
\cite{FF,EF} that as Poisson algebra $Z(\G)_{-\k}$ is isomorphic to the
classical ${\cal W}$--algebra ${\cal W}(^L\!\g)$ associated with the
Langlands dual Lie algebra $\gL$ of $\g$.

In \cite{FR} two of the authors used this second realization of
$\W$--algebras to obtain their $q$--deformations. For instance the
$q$--deformation $\W_q(\sw_n)$ of ${\cal W}(\sw_n)$ was defined as the
center $Z_q(\sun)$ of a completion of the quantized universal
enveloping algebra $U_q(\sun)_{-\k}$. The Poisson structure on
$Z_q(\sun)$ was explicitly described in \cite{FR} using results of
\cite{ext}. It was shown that the underlying Poisson manifold of
$Z_q(\sun)={\cal W}_q(\sw_n)$ is the space of $q$--difference
operators of the form $D_q^n + t_1 D_q^{n-1} + \ldots + t_{n-1} D_q +
1$. Furthermore, in \cite{FR} a $q$--deformation of the Miura
transformation, i.e. a homomorphism from $\W_q(\sw_n)$ to a
Heisenberg-Poisson algebra was defined. The construction \cite{FR} of
${\cal W}_q(\sw_n)$ was followed by further developments: it was
quantized \cite{SKAO,ell,AKOS} and the quantum algebra was used in the
study of lattice models \cite{LP,Miwa}; Yangian analogue of ${\cal
W}_q(\sw_2)$ was considered in \cite{Hou}; $q$--deformations of the
generalized KdV hierarchies were introduced \cite{Fr}.

In this paper we first formulate the results of \cite{FR} in terms of
first order $q$--difference operators and $q$--gauge action. This
naturally leads us to a generalization of the Drinfeld-Sokolov scheme
to the setting of $q$--difference operators. The initial Poisson
manifold is the loop group $LSL_n$ of $SL_n$, or more generally, the
loop group of a simply-connected simple Lie group $G$. Much of the
needed Poisson formalism has already been developed by one of the
authors in \cite{RIMS,dual}. Results of these works allow us to
define a Poisson structure on the loop group, with respect to which the
$q$--gauge action is Poisson. We then have to perform a reduction of
this Poisson manifold with respect to the $q$--gauge action of the
loop group $LN$ of the unipotent subgroup $N$ of $G$.

At this point we encounter a new kind of anomaly in the Poisson
bracket relations, unfamiliar from the linear, i.e. undeformed,
situation. To describe it in physical terms, recall that the reduction
procedure consists of two steps: (1) imposing the constraints and (2)
passing to the quotient by the gauge group. An important point in the
ordinary Drinfeld-Sokolov reduction is that these constraints are of
first class, according to Dirac, i.e., their Poisson bracket vanishes
on the constraint surface. In the $q$--difference case we have to
choose carefully the classical $r$--matrix defining the initial
Poisson structure on the loop group so as to make all constraints
first class. If we use the standard $r$--matrix, some of the
constraints are of second class, and so we have to modify it.

In this paper we do that in the case of $SL_2$. We show that there is
essentially a unique classical $r$--matrix compatible with the
$q$--difference Drinfeld-Sokolov scheme. To the best of our knowledge,
this classical $r$--matrix is new; it yields an elliptic deformation
of the Lie bialgebra structure on the loop algebra of $\sw_2$
associated with the Drinfeld ``new'' realization of quantized affine
algebras \cite{nr}, \cite{kh-t}. The result of the Drinfeld-Sokolov
reduction is the $q$--deformation of the classical Virasoro algebra
defined in \cite{FR}.

We also construct a finite difference version of the Drinfeld-Sokolov
reduction in the case of $SL_2$. This construction gives us a discrete
version of the (classical) Virasoro algebra. We explain in detail the
connection between our discrete Virasoro algebra and the lattice
Virasoro algebra of Faddeev--Takhtajan--Volkov \cite{ft,v1,v2,fv}. We
hope that our results will help to clarify further the meaning of the
discrete Virasoro algebra and its relation to various integrable
models.

The construction presented here can be generalized to the case of an
arbitrary simply-connected simple Lie group. This is done in the
second part of the paper \cite{SS} written by A.~Sevostyanov and one
of us.

The paper is arranged as follows. In Sect.~2 we recall the relevant
facts of \cite{DS} and \cite{FR}. In Sect.~3 we interpret the results
of \cite{FR} from the point of view of $q$--gauge transformations.
Sect.~4 reviews some background material on Poisson structures on Lie
groups following \cite{RIMS,dual}. In Sect.~5 we apply the results of
Sect.~4 to the $q$--deformation of the Drinfeld-Sokolov reduction in
the case of $SL_2$. In Sect.~6 we discuss the finite difference
analogue of this reduction and compare its results with the
Faddeev-Takhtajan-Volkov algebra.

\subsubsection*{Acknowledgements.} E.~Frenkel thanks
P.~Schapira for his hospitality at Universit\'{e} Paris VI in June of
1996, when this collaboration began. Some of the results of this paper
have been reported in E.~Frenkel's lecture course on Soliton Theory
given at Harvard University in the Spring of 1996.

The research of E.~Frenkel was supported by grants from the Packard
and Sloan Foundations, and by the NSF grants DMS 9501414 and DMS
9304580. The research of N.~Reshetikhin was supported by the NSF grant
DMS 9296120.

\section{Preliminaries}

\subsection{The differential Drinfeld-Sokolov reduction in the case of
$\sw_n$}    \label{diff1}

Let ${\cal M}_n$ be the manifold of differential operators of the form
\beq
\label{L} L = \pa^n + u_1(s) \pa^{n-2} + \ldots +
u_{n-2}(s) \pa + u_{n-1}(s),
\end{equation}
where $u_i(s) \in \C((s))$.

Adler \cite{A} and Gelfand-Dickey \cite{GD} have defined a remarkable
two-parameter family of Poisson structures on ${\cal M}_n$, with
respect to which the corresponding KdV hierarchy is hamiltonian. In
this paper we will only consider one of them, the so-called second
bracket. There is a simple realization of this structure in terms of
the Drinfeld-Sokolov reduction \cite{DS}, Sect.~6.5. Let us briefly
recall this realization.

Consider the affine Kac-Moody algebra $\sun$ associated to $\sw_n$;
this is the central extension $$0 \arr \C K \arr \sun \arr \lsln \arr
0,$$ see \cite{Kac}. Let $M_n$ be the hyperplane in the dual space to
$\sun$, which consists of linear functionals taking value $1$ on
$K$. Using the differential $dt$ and the bilinear form $\on{tr} AB$ on
$\sw_n$, we identify $M_n$ with the manifold of first order
differential operators
$$\pa_s + A(s), \quad \quad A(s) \in \lsln.$$ The coadjoint action of
the Lie group $\wh{SL}_n$ on $\sun^*$ factors through the loop group
$LSL_n$ and preserves the hyperplane $M_n$. The corresponding action
of $g(s) \in LSL_n$ on $M_n$ is given by \beq \label{coadjoint} g(s)
\cdot (\pa_s + A(s)) = g(s)(\pa_s + A(s))g(s)^{-1},
\end{equation} or $$A(s) \mapsto g(s) A(s) g(s)^{-1} - \pa_s g(s) \cdot
g(s)^{-1}.$$

Consider now the submanifold $M^J_n$ of $M_n$ wich consists of
operators $\pa_s + A(s)$, where $A(s)$ is a traceless matrix of the form
\beq    \label{upper}
\begin{pmatrix}
* & * & * & \hdots & * & * \\
-1 & * & * & \hdots & * & * \\
0 & -1 & * & \hdots & * & * \\
\hdotsfor{6} \\
0 & 0 & 0 & \hdots & * & * \\
0 & 0 & 0 & \hdots & -1 & *
\end{pmatrix}
\end{equation}

To each element ${\cal L}$ of $M^J_n$ one can naturally attach an
$n$th order scalar differential operator as follows. Consider the
equation ${\cal L} \cdot \Psi = 0$, where $$\Psi = \begin{pmatrix}
\Psi_n \\ \Psi_{n-1} \\ \hdots \\ \Psi_1 \end{pmatrix}.$$
Due to the special form \eqref{upper} of ${\cal L}$, this equation is
equivalent to an $n$th order differential equation $L \cdot \Psi_1 =
0$, where $L$ is of the form \eqref{L}. Thus, we obtain a map $\pi:
M^J_n \arr \M$ sending ${\cal L}$ to $L$.

Let $N$ be the subgroup of $SL_n$ consisting of the upper triangular
matrices, and $LN$ be its loop group. If $g \in LN$ and $\Psi$ is a
solution of ${\cal L} \cdot \Psi = 0$, then $\Psi'= g \Psi$ is a
solution of ${\cal L}' \cdot \Psi' = 0$, where ${\cal L}' = g {\cal L}
g^{-1}$. But $\Psi_1$ does not change under the action of
$LN$. Therefore $\pi({\cal L}') = \pi({\cal L})$, and we see that
$\pi$ factors through the quotient of $M^J_n$ by the action of
$LN$. The following proposition describes this quotient.

\begin{prop}[\cite{DS}, Proposition 3.1]    \label{free}
The action of $LN$ on $M^J_n$ is free, and each orbit contains a
unique operator of the form
\beq    \label{special}
\pa_s + \begin{pmatrix}
0 & u_1 & u_2 & \hdots & u_{n-2} & u_{n-1} \\
-1 & 0 & 0 & \hdots & 0 & 0 \\
0 & -1 & 0 & \hdots & 0 & 0 \\
\hdotsfor{6} \\
0 & 0 & 0 & \hdots & 0 & 0 \\
0 & 0 & 0 & \hdots & -1 & 0
\end{pmatrix}.
\end{equation}
\end{prop}

But for ${\cal L}$ of the form \eqref{special}, $\pi({\cal L})$ is
equal to the operator $L$ given by formula \eqref{L}. Thus, we have
identified the map $\pi$ with the quotient of $M^J_n$ by $LN$ and
identified ${\cal M}_n$ with $M^J_n/LN$.

The quotient $M^J_n/LN$ can actually be interpreted as the result of
hamiltonian reduction. The manifold $M_n$ has a canonical Poisson
structure, which is the restriction of the Lie-Poisson structure on
$\sun^*$ (such a structure exists on the dual space to any Lie
algebra). The coadjoint action of $LN$ on $M_n$ is hamiltonian with
respect to this structure. The corresponding moment map $\mu: M_n \arr
L\n_- \simeq L\n^*$ sends $\pa_s + A(s)$ to the lower-triangular part
of $A(s)$. Consider the one-point orbit of $LN$,
$$J = \begin{pmatrix} 0 & 0 & 0 & \hdots & 0 & 0 \\ -1 & 0 & 0 &
\hdots & 0 & 0 \\ 0 & -1 & 0 & \hdots & 0 & 0 \\ \hdotsfor{6} \\ 0 & 0
& 0 & \hdots & 0 & 0 \\ 0 & 0 & 0 & \hdots & -1 & 0
\end{pmatrix}.$$ Then $M^J_n = \mu^{-1}(J)$. Hence $M_n$ is the
result of hamiltonian reduction of $M_n$ by $LN$ with respect to the
one-point orbit $J$.

The Lie-Poisson structure on $M_n$ gives rise to a canonical Poisson
structure on ${\cal M}_n$, which coincides with the second
Adler-Gelfand-Dickey bracket, see \cite{DS}, Sect.~6.5. The Poisson
algebra of local functionals on ${\cal M}_n$ is called the classical
$\W$--algebra associated to $\sw_n$, and is denoted by $\W(\sw_n)$.

\begin{rem}    \label{kostant}
For $\al \in \C$, let $M_{\al,n}$ be the hyperplane in the dual space
to $\sun$, which consists of linear functionals on $\sun$ taking value
$\al$ on $K$. In the same way as above (for $\al=1$) we identify
$M_{\al,n}$ with the space of first order differential operators
$$\al \pa_s + A(s), \quad \quad A(s) \in \lsln.$$ The coadjoint action
is given by the formula $$A(s) \mapsto g(s) A(s) g(s)^{-1} - \al \pa_s
g(s) \cdot g(s)^{-1}.$$ The straightforward generalization of
\propref{free} is true for any $\al \in \C$. In particular, for
$\al=0$ we obtain a description of the orbits in $M_n^J$ under the
adjoint action of $LN$. This result is due to B.~Kostant \cite{Ko}.

Drinfeld and Sokolov \cite{DS} gave a generalization of \propref{free}
when $\sw_n$ is replaced by an arbitrary semisimple Lie algebra
$\g$. The special case of their result, corresponding to $\al=0$,
is also due to Kostant \cite{Ko}.\qed
\end{rem}

The Drinfeld-Sokolov reduction can be summarized by the following
diagram.

\setlength{\unitlength}{1mm}

\begin{center}
\begin{picture}(125,40)(-35,-33)     
\put(0,-1){\vector(0,-1){26}}
\put(62,-1){\vector(0,-1){26}}
\put(9,3){\vector(1,0){45}}
\put(9,-31){\vector(1,0){45}}
\put(8.9,3.9){\oval(1,1)[l]}
\put(8.9,-30.1){\oval(1,1)[l]}
\put(-2,2){$M_n^J$}
\put(-23,-32){$\M=M_n^J/LN$}
\put(59,2){$M_n$}
\put(57,-32){$M_n/LN$}

\end{picture}
\end{center}

There are three essential properties of the Lie-Poisson structure on
$M_n$ that make the reduction work:

\begin{itemize}

\item[(i)] The coadjoint action of $LSL_n$ on $M_n$ is hamiltonian
with respect to this structure;

\item[(ii)] the subgroup $LN$ of $LSL_n$ is admissible in the sense that
the space $S$ of $LN$--invariant functionals on $M_n$ is a Poisson
subalgebra of the space of all functionals on $M_n$;

\item[(iii)] Denote by $\mu_{ij}$ the function on $M_n$, whose
value at $\pa + A \in M_n$ equals the $(i,j)$ entry of $A$. The ideal
in $S$ generated by $\mu_{ij} + \delta_{i-1,j}, i>j$, is a Poisson
ideal.

\end{itemize}

We will generalize this picture to the $q$--difference case.

\subsection{The Miura transformation}

Let $\F$ be the manifold of differential operators of the
form
\begin{equation}    \label{cartan}
\pa_s +
\begin{pmatrix}
v_1 & 0 & 0 & \hdots & 0 & 0 \\
-1 & v_2 & 0 & \hdots & 0 & 0 \\
0 & -1 & v_3 & \hdots & 0 & 0 \\
\hdotsfor{6} \\
0 & 0 & 0 & \hdots & v_{n-1} & 0 \\
0 & 0 & 0 & \hdots & -1 & v_n
\end{pmatrix},
\end{equation}
where $\sum_{i=1}^n v_i = 0$.

We have a map $\muu: \F \arr \M$, which is the composition of the
embedding $\F \arr M^J_n$ and the projection $\pi: M^J_n \arr \M$.

Using the definition of $\pi$ above, $\muu$ can be described explicitly
as follows: the image of the operator \eqref{cartan} under $\muu$ is
the $n$th order differential operator
$$\pa_s^n + u_1(s) \pa_s^{n-2} + \ldots + u_{n-1}(s) = (\pa_s +
v_1(s)) \ldots (\pa_s + v_n(s)).$$

The map $\muu$ is called the Miura transformation.

We want to describe the Poisson structure on $\F$ with respect to
which the Miura transformation is Poisson. To this end, let us
consider the restriction of the gauge action (\ref{coadjoint}) to the
opposite triangular subgroup $LN_{-};$ let $\ol{\mu} :M_n\rightarrow
L {\frak n}_+ \simeq L{\frak n}_{-}^{*}$ be the corresponding moment
map. The manifold $\F$ is the intersection of two level surfaces,
$\F = \mu^{-1}(J)\cap \ol{\mu}^{-1}(0).$ It is easy to see that it
gives a local cross-section for both actions (in other words, the
orbits of $LN$ and $LN_{-}$ are transversal to $F_n$). Hence $F_n$
simultaneously provides a local model for the reduced spaces
$M_n=\mu^{-1}(J)/LN$ and $\ol{\mu}^{-1}(0)/LN_{-}$. The Poisson
bracket on $\F$ that we need to define is the so-called Dirac
bracket (see, e.g., \cite{Flato}), where we may regard the matrix
coefficients of $\ol{\mu}$ as subsidiary conditions, which fix the
local gauge. The computation of the Dirac bracket for the diagonal
matrix coefficients $v_i$ is very simple, since their Poisson brackets
with the matrix coefficients of $\ol{\mu}$ all vanish on $\F$. The
only correction arises due to the constraint $\sum_{i=1}^N v_i = 0$.

Denote by $v_{i,n}$ the linear functional on $\F$, whose value on the
operator \eqref{cartan} is the $n$th Fourier coefficient of
$v_i(s)$. We obtain the following formula for the Dirac bracket on
$\F$:
\begin{align*}
\{ v_{i,n},v_{i,m} \} &= \frac{N-1}{N} n \delta_{n,-m}, \\ \{
v_{i,n},v_{j,m} \} &= - \frac{1}{N} n \delta_{n,-m}, \quad i<j.
\end{align*}

Since $\F$ and $\M$ both are models of the same reduced space, we
immediately obtain:

\begin{prop}[\cite{DS}, Proposition 3.26]
With respect to this Poisson structure the map $\muu: \F \arr \M$ is
Poisson.
\end{prop}

\subsection{The $q$--deformations of $\W(\sw_n)$ and Miura transformation}
\label{qdef}

In this section we summarize relevant results of \cite{FR}.

Consider the space ${\cal M}_{n,q}$ of $q$--difference operators of
the form \beq \label{qL} L = D^n + t_1(s) D^{n-1} + \ldots +
t_{n-1}(s) D + 1,
\end{equation}
where $t_i(s) \in \C((s))$ for each $i=1,\ldots,n$, and $[D \cdot
f](s) = f(sq)$.

Denote by $t_{i,m}$ the functional on $\Mq$, whose value at $L$ is the
$m$th Fourier coefficient of $t_i(s)$. Let $\RN$ be the completion
of the ring of polynomials in $t_{i,m}, i=1,\ldots,N-1; m \in \Z$,
which consists of finite linear combinations of expressions of the
form \beq
\label{qloc} \sum_{m_1+\ldots +m_k=M} c(m_1,\ldots,m_k) \cdot
t_{i_1,m_1} \ldots t_{i_k,m_k},
\end{equation}
where $c(m_1,\ldots,m_k) \in \C$. Given an operator of the form
\eqref{qL}, we can substitute the coefficients $t_{i,m}$ into an
expression like \eqref{qloc} and get a number. Therefore elements of
$\RN$ define functionals on the space $\Mq$.

In order to define the Poisson structure on $\Mq$, it suffices to
specify the Poisson brackets between the generators $t_{i,m}$. Let
$T_i(z)$ be the generating series of the functionals $t_{i,m}$:
$$T_i(z) = \sum_{m\in\Z} t_{i,m} z^{-m}.$$ We define the Poisson
brackets between $t_{i,m}$'s by the formulas \cite{FR}
\begin{align} \notag
\{ T_i(z),T_j(w) \} = & \sum_{m\in\Z} \left( \frac{w}{z}
\right)^m \frac{(1-q^{im})(1-q^{m(N-j)})}{1-q^{mN}} T_i(z) T_j(w) \\
\label{p2} &+ \sum_{r=1}^{\on{min}(i,N-j)} \delta \left( \frac{wq^r}{z}
\right) T_{i-r}(w) T_{j+r}(z) \\ \notag &- \sum_{r=1}^{\on{min}(i,N-j)}
\delta\left( \frac{w}{zq^{j-i+r}} \right) T_{i-r}(z) T_{j+r}(w), \quad
i\leq j.
\end{align}
In these formulas $\delta(x) = \sum_{m\in\Z} x^m$, and we use the
convention that $t_0(z) \equiv 1$.

\begin{rem}
The parameter $q$ in formula \eqref{p2} corresponds to $q^{-2}$ in the
notation of \cite{FR}.\qed
\end{rem} \medskip

Now consider the space $\Fq$ of $n$--tuples of $q$--difference operators
\beq    \label{ntuple}
(D+\la_1(s),\ldots,D+\la_n(s)).
\end{equation}
Denote by $\la_{i,m}$ the functional on $\Mq$, whose value is the
$m$th Fourier coefficient of $\la_i(s)$. We will denote by
$\La_i(z)$ the generating series of the functionals $\la_{i,m}$:
$$\La_i(z) = \sum_{m\in\Z} \la_{i,m} z^{-m}.$$

We define a Poisson structure on $\Fq$ by the formulas \cite{FR}:
\begin{align}    \label{pbg1}
\{ \La_i(z),\La_i(w) \} & = \sum_{m\in\Z} \left( \frac{w}{z} \right)^m
\frac{(1-q^m)(1-q^{m(N-1)})}{1-q^{mN}} \La_i(z) \La_i(w), \\ \label{pbg2}
\{ \La_i(z),\La_j(w) \} & = - \sum_{m\in\Z} \left(
\frac{wq^{N-1}}{z} \right)^m \frac{(1-q^m)^2}{1-q^{mN}} \La_i(z)
\La_j(w), \quad i<j.
\end{align}

Now we define the $q$--deformation of the Miura transformation as the
map $\muu_q: \Fq \arr \Mq$, which sends the $n$--tuple \eqref{ntuple}
to \beq
\label{Li} L = (D+\la_1(s))(D+\la_2(sq^{-1})) \ldots
(D+\la_N(sq^{-N+1})),
\end{equation}
i.e.
\beq    \label{formulai}
t_i(s) = \sum_{j_1 < \ldots < j_i} \la_{j_1}(s) \la_{j_2}(sq^{-1})
\ldots \la_{j_i}(sq^{-i+1}).
\end{equation}

\begin{prop}[\cite{FR}]
The map $\muu_q$ is Poisson.
\end{prop}

\subsection{$q$--deformation of the Virasoro algebra}
Here we specialize the formulas of the previous subsection to the
case of $\sw_2$ (we will omit the index $1$ in these formulas). We
have the following Poisson bracket on $T(z)$:
\beq    \label{vira}
\{ T(z),T(w) \} = \sum_{m\in\Z} \left( \frac{w}{z}
\right)^m \frac{1-q^m}{1+q^m} T(z) T(w) + \delta \left( \frac{wq}{z}
\right) - \delta \left( \frac{w}{zq} \right).
\end{equation}

The $q$--deformed Miura transformation reads:
\beq    \label{virmiura}
\La(z) \mapsto T(z) = \La(z) + \La(zq)^{-1}.
\end{equation}

The Poisson bracket on $\La(z)$:
\beq    \label{virpois}
\{ \La(z),\La(w) \} = \sum_{m\in\Z} \left( \frac{w}{z}
\right)^m \frac{1-q^m}{1+q^m} \La(z) \La(w).
\end{equation}

\section{Connection with $q$--gauge transformations}    \label{qgauge}

In this section we present the results of \cite{FR} in the form
of $q$--difference Drinfeld-Sokolov reduction.

\subsection{Presentation via first order $q$--difference operators}

By analogy with the differential case, it is natural to consider the
manifold $M_{n,q}$ of first order difference operators $D + A(s)$,
where $A(s)$ is an element of the loop group $LSL_n$ of $SL_n$. The
group $LSL_n$ acts on this manifold by the $q$--gauge transformations
\beq \label{qadjoint} g(s) \cdot (D + A(s)) = g(sq)(D +
A(s))g(s)^{-1},
\end{equation}
i.e. $g(s) \cdot A(s) = g(sq) A(s) g(s)^{-1}$.

Now we consider the submanifold $M^J_{n,q} \subset M_{n,q}$ which
consists of operators $D+A(s)$, where $A(s)$ is of the form
\eqref{upper}. It is preserved under the $q$--gauge action of the
group $LN$.

In the same way as in the differential case, we define a map $\pi_q:
M^J_{n,q} \arr \Mq$, which sends each element of $M^J_{n,q}$ to an
$n$th order $q$--difference operator $L$ of the form \eqref{qL}. It is
clear that the map $\pi_q$ factors through the quotient of $M^J_{n,q}$
by $LN$. Now we state the $q$--difference analogue of \propref{free}.

\begin{lem}    \label{qfree}
The action of $LN$ on $M^J_{n,q}$ is free and each orbit contains a
unique operator of the form \beq \label{qcan} D +
\begin{pmatrix} t_1 & t_2 & t_3 & \hdots & t_{n-1} & 1 \\ -1 & 0 & 0 &
\hdots & 0 & 0 \\ 0 & -1 & 0 & \hdots & 0 & 0 \\ \hdotsfor{6} \\ 0 & 0
& 0 & \hdots & 0 & 0 \\ 0 & 0 & 0 & \hdots & -1 & 0
\end{pmatrix}.
\end{equation}
\end{lem}

\noindent {\em Proof} is an exercise in elementary matrix algebra. For
$\al=1,\ldots,n$, denote by $M^\al_{n,q}$ the subset of matrices from
$M^J_{n,q}$ satisfying the property that all entries in their rows $i=
\al+1,\ldots,n$ are zero except for the $(i,i-1)$ entry equal to
$-1$. We will prove that given $A(s) \in M^\al_{n,q}, \al>1$, there
exists $g(s) \in LN$, such that $g(sq) A(s) g(s)^{-1} \in
M^{\al-1}_{n,q}$. Since the condition is vacuous for $\al=n$, i.e.
$M^n_{n,q} = M^J_{n,q}$, this will imply that each $LN$--orbit in
$M^J_{n,q}$ contains an element of the form \eqref{qcan}.

To prove the statement for a given $\al$, we will recursively
eliminate all entries of the $\al$th row of $A(s)$ (except the
$(\al,\al-1)$ entry), from right to left using elementary unipotent
matrices. Denote by $E_{i,j}(x)$ the upper unipotent matrix whose only
non-zero entry above the diagonal is the $(i,j)$ entry equal to
$x$. At the first step, we eliminate the $(\al,n)$ entry $A_{\al,n}$
of $A(s)$ by applying the $q$--gauge transformation \eqref{qadjoint}
with $g(s) = E_{\al-1,n}(-A_{\al,n}(s))$. Then we obtain a new matrix
$A'(s)$, which still belongs to $M^\al_{n,q}$, but whose $(\al,n)$
entry is equal to $0$. Next, we apply the $q$--gauge transformation by
$E_{\al-1,n-1}(-A'_{\al,n-1}(s))$ to eliminate the $(\al,n-1)$ entry
of $A'(s)$, etc. It is clear that at each step we do not spoil the
entries that have already been set to $0$. The product of the
elementary unipotent matrices constructed at each step gives us an
element $g(s) \in LN$ with the desired property that $g(sq) A(s)
g(s)^{-1} \in M^{\al-1}_{n,q}$.

To complete the proof, it suffices to remark that if $A(s)$ and
$A'(s)$ are of the form \eqref{qcan}, and $g(sq) A(s) g(s)^{-1} =
A'(s)$ for some $g(s) \in LN$, then $A(s)=A'(s)$ and $g(s)=1$.\qed

For ${\cal L}$ of the form \eqref{special}, $p({\cal L})$ equals the
operator $L$ given by formula \eqref{L}. Thus, we have identified the
map $\pi_q$ with the quotient of $M^J_{n,q}$ by $LN$ and $\Mq$ with
$M^J_{n,q}/LN$.

\begin{rem}
In the same way as above we can prove the following more general
statement. Let $R$ be a ring with an automorphism $\tau$. It rives
rise to an automorphism of $SL_n(R)$ denoted by the same
character. Define $M^J_{\tau,n}$ as the set of elements of $SL_n(R)$
of the form \eqref{upper}. Let the group $N(R)$ act on
$M^J_{\tau,n}(R)$ by the formula $g \cdot A = (\tau \cdot g) A
g^{-1}$. Then this action of $N(R)$ is free, and the quotient is
isomorphic to the set ${\cal M}^J_{\tau,n}(R)$ of elements of
$SL_n(R)$ of the form \eqref{qcan} (i.e. each orbit contains a unique
element of the form \eqref{qcan}). Note that the proof is not sensible
to whether $\tau=\on{Id}$ or not.

When $\tau=\on{Id}$, this result is well-known. It gives the classical
normal form of a linear operator. Moreover, in that case R.~Steinberg
has proved that the subset ${\cal M}^J_{\on{Id},n}(K)$ of $SL_n(K)$,
where $K$ is an algebraically closed field, is a cross-section of the
collection of regular conjugacy classes in $SL_n(K)$ \cite{St},
Theorem 1.4. Steinberg defined an analogous cross-section for any
simply-connected semisimple algebraic group \cite{St}. His results can
be viewed as group analogues of Kostant's results on semisimple Lie
algebras \cite{Ko} (cf. \remref{kostant}). Steinberg's cross-section
is used in the definition of the discrete Drinfeld-Sokolov reduction
in the general semisimple case (see \cite{SS}).\footnote{We are
indebted to B.~Kostant for drawing our attention to \cite{St}}\qed
\end{rem}

\subsection{Deformed Miura transformation via $q$--gauge action}

Let us attach to each element of $\Fq$ the $q$--difference operator
\beq    \label{Lambda}
\La = D + \begin{pmatrix}
\la_1(s) & 0 & \hdots & 0 & 0 \\
-1 & \la_2(sq^{-1}) & \hdots & 0 & 0 \\
\hdotsfor{5} \\
0 & 0 & \hdots & \la_{n-1}(sq^{-n+2}) & 0 \\
0 & 0 & \hdots & -1 & \la_n(sq^{-n+1})
\end{pmatrix},
\end{equation}
where $\prod_{i=1}^n \la_i(sq^{-i+1}) = 1$.

Let $\wt{\muu}_q: \Fq \arr \Mq$ be the composition of the embedding
$\Fq \arr M^J_{n,q}$ and $\pi_q: M^J_{n,q} \arr M^J_{n,q}/LN \simeq
\Mq$. Using the definition of $\pi_q$ above, one easily finds that for
$\La$ given by \eqref{Lambda}, $\wt{\muu}_q(\La)$ is the operator
\eqref{qcan}, where $t_i(s)$ is given by formula \eqref{formulai}.

Therefore we obtain

\begin{lem}
The map $\wt{\muu}_q$ coincides with the $q$--deformed Miura
transformation $\muu_q$.
\end{lem}

\begin{rem}    \label{qchar}
Let $G$ be a simply-connected semisimple algebraic group over $\C$.
Let $V_i$ be the $i$th fundamental representation of $G$ (in the case
$G=SL_n$, $V_i = \Lambda^i \C^n$), and $\chi_i: G \arr \C$ be the
corresponding character, $\chi_i(g) = \on{Tr}(g,V_i)$. Define a map
$p: G \arr \C^N$ by the formula $p(g) =
(\chi_1(g),\ldots,\chi_n(g))$. By construction, $p$ is constant on
conjugacy classes. In the case $G=SL_n$ the map $p$ has a
cross-section $r: \C^n \arr SL_n(\C)$:
$$(a_1,\ldots,a_n) \mapsto \begin{pmatrix} a_1 & a_2 & a_3 & \hdots &
a_{n-1} & 1 \\ -1 & 0 & 0 & \hdots & 0 & 0 \\ 0 & -1 & 0 & \hdots & 0
& 0 \\ \hdotsfor{6} \\ 0 & 0 & 0 & \hdots & 0 & 0 \\ 0 & 0 & 0 &
\hdots & -1 & 0
\end{pmatrix}.$$ The composition $r \circ p$, restricted to
$M^J_{n,1}$ coincides with the map $\pi_1$. Moreover, $\wt{\muu}_1$
can be interpreted as the restriction of $p$ to the subset of $SL_n$
consisting of matrices of the form \eqref{Lambda}. Hence $\wt{\muu}_1$
sends $(\la_1,\ldots,\la_n)$ to the elementary symmetric polynomials
$$t_i = \sum_{j_1 < \ldots < j_i} \la_{j_1} \la_{j_2} \ldots
\la_{j_i},$$ which are the characters of the fundamental
representations of $SL_n$. As we mentioned above, Steinberg has defined
an analogue of the cross-section $r$ for an arbitrary simply-connected
semisimple algebraic group \cite{St}.

Formula \eqref{formulai} means that in terms of $\la_j(z)$ the
generators $t_i(z)$ of $W_q(\sw_n)$ can be thought of as
$q$--deformations of the characters of fundamental representations of
$SL_n$. It is interesting that the same interpretation is also
suggested by the definintion of $\W_q(\sw_n)$ as the center of a
completion of the quantized universal enveloping algebra
$U_q(\sun)_{-\k}$ \cite{FR}. Namely, $t_i(z)$ is then defined as the
($q$--deformed) trace of the so-called $L$--operator acting on
$\Lambda^i \C^n$ considered as a representation of $U_q(\sun)$, see
\cite{ext,FR} (note also that $t_i(z)$ is closely connected with a
transfer-matrix of the corresponding integrable spin model).\qed
\end{rem}

Thus, we have now represented $\Mq$ as the quotient of the submanifold
$M^J_{n,q}$ of the manifold $M_{n,q}$ of first order $q$--difference
operators by the action of the group $LN$ (acting by $q$--gauge
transformations). We have also interpreted the $q$--deformed Miura
transformation in these terms. In the next sections we discuss the
Poisson structure on $M_{n,q}$, which gives rise to the Poisson
structure on $\Mq$ given by explicit formula \eqref{p2}.

\section{Poisson structures}    \label{P}

\subsection{Overview}    \label{over}
In view of the previous section, the following diagram is the
$q$--difference analogue of the diagram presented at the end of
\secref{diff1}.

\begin{center}
\begin{picture}(125,40)(-30,-33)     
\put(0,-1){\vector(0,-1){26}}
\put(62,-1){\vector(0,-1){26}}
\put(9,3){\vector(1,0){45}}
\put(9,-31){\vector(1,0){45}}
\put(8.9,3.9){\oval(1,1)[l]}
\put(8.9,-30.1){\oval(1,1)[l]}
\put(-5,2){$M_{n,q}^J$}
\put(-26,-32){$\Mq=M_{n,q}^J/LN$}
\put(57,2){$M_{n,q}=(LG,\,\eta_*^q)$}
\put(57,-32){$M_{n,q}/LN$}

\end{picture}
\end{center}

As in the differential case, in order to define a $q$--deformation of
the Drinfeld-Sokolov reduction we need to find a Poisson structure
$\eta_*^q$ on $M_{n,q}$ and a Poisson-Lie structure $\eta$ on $LSL_n$
satisfying the following properties.

\begin{itemize}

\item[(i)] the action $LSL_n \times M_{n,q} \arr M_{n,q}$ by
$q$--gauge transformations is Poisson;

\item[(ii)] the subgroup $LN$ of $LSL_n$ is admissible in the sense
that the algebra $S_q$ of $LN$--invariant functionals on $M_{n,q}$
is a Poisson subalgebra of the algebra of all functionals on $M_{n,q}$;

\item[(iii)] Denote by $\mu_{ij}$ the function on $M_{n,q}$, whose
value at $D + A \in M_{n,q}$ equals the $(i,j)$ entry of $A$. The
ideal in $S_q$ generated by $\mu_{ij} + \delta_{i-1,j}, i>j$, is a
Poisson ideal.

\end{itemize}

Geometrically, the last condition means that $\Mq$ is a Poisson
submanifold of the quotient $M_{n,q}/LN$.

For the sake of completeness, we recall the notions mentioned
above. Let $M$ be a Poisson manifold, and $H$ be a Lie group, which is
itself a Poisson manifold. An action of $H$ on $M$ is called Poisson
if $H \times M\rightarrow M$ is a Poisson map (here we equip $H\times
M$ with the product Poisson structure). In particular, if the
multiplication map $H \times H \arr H$ is Poisson, then $H$ is called
a Poisson-Lie group.

In this section we describe the general formalism concerning problems
(i)--(iii) above. Then in the next section we specialize to $M_{2,q}$
and give an explicit solution of these problems.

\subsection{Lie bialgebras}    \label{bialg}
Let $\g$ be a Lie algebra. Recall \cite{Dr} that $\g$ is called a Lie
bialgebra, if $\g^*$ also has a Lie algebra structure, such that the dual
map $\delta: \g \arr \Lambda^2 \g$ is a one-cocycle. We will consider {\em
factorizable} Lie bialgebras ($\g,\delta$) satisfying the following
conditions:

\begin{itemize}

\item[(1)] There exists a linear map $r_+: \g^* \arr \g$, such that both
$r_+$ and $r_-= -r_+^*$ are Lie algebra homomorphisms.

\item[(2)] The endomorphism $t = r_+ -r_-$ is $\g$-equivariant and
induces a linear isomorphism $\g^*\arr\g$.

\end{itemize}

Instead of the linear operator $r_+ \in \on{Hom}(\g^*,\g)$ one often
considers the corresponding element $r$ of $\g^{\ot 2}$ (or a
completion of $\g^{\ot 2}$ if $\g$ is infinite-dimensional). The
element $r$ (or its image in the tensor square of a particular
representation of $\g$) is called classical $r$--matrix. It satisfies
the classical Yang-Baxter equation:
\begin{equation}    \label{yb}
[r_{12},r_{13}] + [r_{12},r_{23}] + [r_{13},r_{23}] = 0.
\end{equation}
In terms or $r$, $\delta(x)=[r,x], \forall x \in \g$ (here $[a \otimes
b,x] = [a,x] \otimes b + a \otimes [b,x]$). The maps $r_\pm: \g^* \arr
\g$ are given by the formulas: $r_+(y) = (y \otimes \on{id})(r),
r_-(y) = - (\on{id} \otimes y)(r)$.

Property (2) above means that $r+\sigma(r)$, where $\sigma(a \otimes
b) = b \otimes a$ is a non-degenerate $\g$--invariant symmetric
bilinear form on $\g^*$.

Set $\g_\pm = \on{Im}(r_\pm)$. Property (1) above implies that
$\g_\pm \subset \g$ is a Lie subalgebra. The following statement is
essentially contained in \cite{BD} (cf. also \cite{rmatr}).

\begin{lem}
Let $(\g, \g^*)$ be a factorizable Lie bialgebra. Then

(1) The subspace $\n_\pm = r_\pm(\on{Ker} r_\mp)$ is a Lie ideal in
$\g_\pm$.

(2) The map $\theta: \g_+/\n_+\arr\g_-/\n_-$ which sends the
residue class of $r_+(X), X\in\g^*$, modulo $\n_+$ to that of $r_-(X)$
modulo $\n_-$ is a well-defined isomorphism of Lie algebras.

(3) Let $\D=\g\oplus\g$ be the direct sum of two copies of
$\g$. The map $$i: \g^* \arr \D, \quad \quad X \mapsto
(r_+(X),r_-(X))$$ is a Lie algebra embedding; its image $\g^* \subset
\D$ is $$\g^* = \{(X_+,X_-) \in \g_+ \oplus \g_- \subset \D | \ol{X}_- =
\theta(\ol{X}_+) \},$$ where $\ol{Y}_\pm = Y \on{mod} \n_\pm$.
\end{lem}

\begin{rem} The connection between our notation and that of
\cite{RIMS} is as follows: the operator $r \in \on{End}\g$ of
\cite{RIMS} coincides with the composition of $r_+ + r_-$ up to the
isomorphism $t = r_+ - r_-: \g^* \arr \g$; the bilinear form used in
\cite{RIMS} is induced by $t$.\qed
\end{rem}

\subsection{Poisson-Lie groups and gauge transformations}
\label{adj}
Let ($G,\eta$) (resp., \newline ($G^*,\eta^*$)) be a Poisson-Lie group
with factorizable tangent Lie bialgebra ($\g,\delta$) (resp.,
($\g^*,\delta^*$)). Let $G_\pm$ and $N_\pm$ be the Lie subgroups of
$G$ corresponding to the Lie subalgebras $\g_\pm$ and $\n_\pm$. We
denote by the same symbol $\theta$ the isomorphism $G_+/N_+ \arr
G_-/N_-$ induced by $\theta: \g_+/\n_+ \arr \g_-/\n_-$. Then the group
$G^*$ is isomorphic to $$\{ (g_+,g_-) \in G_+ \times G_- |
\theta(\ol{g}_+) = \ol{g}_- \},$$ and we have a map $i: G^* \arr G$
given by $i((g_+,g_-)) = g_+ (g_-)^{-1}$.

Explicitly, Poisson bracket on ($G,\eta$) can be written as follows:
\beq \label{skl} \{ \vf,\psi \} = \langle r,\nb \vf \wedge \nb \psi -
\nb' \vf \wedge \nb' \psi \rangle,
\end{equation}
where for $x \in G$, $\nb \vf(x), \nb' \vf(x) \in \g^*$ are defined by
the formulas:
\begin{align}    \label{nabla}
\langle \nb \vf(x),\xi \rangle &= \frac{d}{dt} \vf\left( e^{t\xi} x
\right)|_{t=0},\\ \langle \nb' \vf(x),\xi \rangle &= \frac{d}{dt}
\vf\left( x e^{t\xi} \right)|_{t=0},
\end{align}
for all $\xi \in \g$. Analogous formula can be written for the
Poisson bracket on ($G^*,\eta^*$). In formula \eqref{skl} we use the
standard notation $a \wedge b = (a \ot b - b \ot a)/2$.

By definition, the action of $G$ on itself by left translations is a
Poisson group action. There is another Poisson structure $\eta_*$ on
$G$ which is covariant with respect to the adjoint action of $G$ on
itself and such that the map $i: (G^*,\eta^*) \arr (G,\eta_*)$ is
Poisson. It is given by the formula \beq \label{another} \{ \vf,\psi
\} = \langle r,\nb \vf \wedge \nb \psi + \nb' \vf \wedge \nb' \psi
\rangle - \langle r, \nb' \vf \ot \nb \psi - \nb' \psi \ot \nb \vf
\rangle.
\end{equation}

\begin{prop}    \label{embed}
(1) The map $i: G^* \arr G$ is a Poisson map between the Poisson
manifolds ($G^*,\eta^*$) and ($G,\eta_*$);

(2) The Poisson structure $\eta_*$ on $G$ is covariant with respect to
the adjoint action, i.e. the map
$$(G,\eta) \times (G,\eta_*) \arr (G,\eta_*): (g,h) \mapsto g h
g^{-1}$$ is a Poisson map.
\end{prop}

These results are proved in \cite{RIMS}, \S~3 (see also \cite{dual},
\S~2), using the notion of the Heisenberg double of $G$. Formula
\eqref{another} can also be obtained directly from the explicit
formulas for the Poisson structure $\eta^*$ and for the embedding $i$.

More generally, let $\tau$ be an automorphism of $G$, such that the
corresponding automorphism of $\g$ satisfies $(\tau \ot \tau)(r) =
r$. Define a twisted Poisson structure $\eta_*^\tau$ on $G$ by the
formula
\begin{align}
\label{mainb1} \{ \vf,\psi \} &= \langle r,\nb \vf \wedge \nb \psi + \nb'
\vf \wedge \nb' \psi \rangle \\ \notag &- \langle (\tau \otimes
\on{id})(r), \nb' \vf \ot \nb \psi - \nb' \psi \ot \nb \vf \rangle,
\end{align}
and the twisted adjoint action of $G$ on itself by the formula $g
\cdot h = \tau(g) h g^{-1}$.

\begin{thm}    \label{act1}
{\em The Poisson structure $\eta_*^\tau$ on $G$ is covariant with
respect to the twisted adjoint action, i.e. the map
$$(G,\eta) \times (G,\eta_*^\tau) \arr (G,\eta_*^\tau): (g,h) \mapsto
\tau(g) h g^{-1}$$ is a Poisson map.}
\end{thm}

This result was proved in \cite{RIMS}, \S~3 (see also \cite{dual},
\S~2), using the notion of the twisted Heisenberg double of $G$. We
will use \thmref{act1} in two cases. In the first, $G$ is the loop
group of a finite-dimensional simple Lie group $\ol{G}$, and $\tau$ is
the automorphism $g(s) \arr g(sq), q \in \C^\times$. In the second, $G
= \ol{G}^{\Z/N\Z}$, and $\tau$ is the automorphism $(\tau(g))_i \arr
g_{i+1}$. In the first case twisted conjugations coincide with
$q$--gauge transformations, and in the second case they coincide with
lattice gauge transformations.

\subsection{Admissibility and constraints}    \label{admcon}

Let $M$ be a Poisson manifold, $G$ a Poisson Lie group and $G \times M
\arr M$ be a Poisson action. A subgroup $H\subset G$ is called
admissible if the space $C^\infty(M)^H$ of $H$--invariant functions on
$M$ is a Poisson subalgebra in the space $C^\infty(M)$ of all
functions on $M$.

\begin{prop}[\cite{RIMS},Theorem 6]    \label{admiss}
Let $\left( {\frak g},{\frak g}^{*}\right) $ be the tangent
Lie bialgebra of $G.$ A connected Lie subgroup $H\subset G$ with Lie algebra 
${\frak h}\subset {\frak g}$ is admissible if ${\frak h}^{\perp }\subset
{\frak g}^{*}$ is a Lie subalgebra.
\end{prop}

In particular, $G$ itself is admissible. Note that $H\subset G$ is a
Poisson subgroup if and only if ${\frak h}^{\perp }\subset {\frak
g}^{*}$ is an ideal; in that case the tangent Lie bialgebra of $H$ is
$\left( {\frak h},{\frak g}^{*}/{\frak h}^{\bot }\right)$.

Let $H\subset G$ be an admissible subgroup, and $I$ be a Poisson ideal
in $C^\infty(M)^H$, i.e. $I$ is an ideal in the ring $C^\infty(M)^H$,
and $\{ f,g \} \in C^\infty(M)^H$ for all $f \in I, g \in
C^\infty(M)^H$. Then $C^\infty(M)^H/I$ is a Poisson algebra.

More geometrically, the Poisson structure on $C^\infty(M)^H/I$ can be
described as follows. Assume that the quotient $M/H$ exists as a smooth
manifold. Then there exists a Poisson structure on $M/H$ such that the
canonical projection $\pi: M\rightarrow M/H$ is a Poisson map. Hamiltonian
vector fields $\xi_\varphi ,\varphi \in \pi ^{*}C^\infty(M/H),$ generate an
integrable distribution ${\frak H}_\pi$ in $TM$. The following result is
straightforward.

\begin{lem}
\label{poisson}
Let $V\subset M$ be a submanifold preserved by $H$. Then $V/H$ is a
Poisson submanifold of $M/H$ if and only if $V$ is an integral
manifold of ${\frak H}_\pi.$
\end{lem}

The integrality condition means precisely that the ideal $I$ of all
$H$--invariant functions on $M$ vanishing on $V$ is a Poisson ideal in
$C^\infty(M)^H$, and $C^\infty(V/H) = C^\infty(V)^H$ $=
C^\infty(M)^H/I$. If this property holds, we will say that the Poisson
structure on $M/H$ can be restricted to $V/H$. 

\begin{center}
\begin{picture}(125,25)(-30,-18)     
\put(0,0){\vector(0,-1){12}}
\put(50,0){\vector(0,-1){12}}
\put(9,3){\vector(1,0){35}}
\put(9,-16){\vector(1,0){35}}
\put(8.9,3.9){\oval(1,1)[l]}
\put(8.9,-15.1){\oval(1,1)[l]}
\put(-1,2){$V$}
\put(-3,-17){$V/H$}
\put(48,2){$M$}
\put(47,-17){$M/H$}

\end{picture}
\end{center}

The Poisson structure on $V/H$ can be described as follows. Let
$N_V\subset T^{*}M\mid_V$ be the conormal bundle of $V$. Clearly,
$T^{*}V\simeq T^{*}M\mid _V/N_V.$ Let $\varphi ,\psi \in C(V)^H$ and
$\overline{d\varphi},\overline{d\psi}\in T^{*}M\mid _V$ be any
representatives of $d\varphi,d\psi \in T^{*}V.$ Let $P_M\in
\bigwedge^2T\,M$ be the Poisson tensor on $M$.

\begin{lem}
\label{Reduce}
We have 
\begin{equation}
\left\{ \varphi ,\psi \right\} = \left\langle
P_M,\overline{d\varphi }\ot \overline{d\psi}\right\rangle ;
\label{Pbr}
\end{equation}
in particular, the right hand side does not depend on the choice of
$\overline{d\varphi},\overline{d\psi}.$
\end{lem}

\begin{rem}
In the case of Hamiltonian action (i.e. when the Poisson structure on
$H$ is trivial), one can construct submanifolds $V$ satisfying the
condition of \lemref{poisson} using the moment map. Although a similar
notion of the nonabelian moment map in the context of Poisson group
theory is also available \cite{Lu}, it is less convenient. The reason
is that the nonabelian moment map is ``less functorial'' than the
ordinary moment map. Namely, if $G\times M\rightarrow M$ is a
Hamiltonian action with moment map $\mu_G: M\rightarrow {\frak
g}^{*},$ its restriction to a subgroup $H\subset G$ is also
Hamiltonian with moment $\mu_H=p \circ \mu_G$ (here $p: {\frak
g}^{*}\rightarrow {\frak h}^{*}$ is the canonical projection). If $G$
is a Poisson-Lie group, $G^{*}$ its dual, $G\times M\rightarrow M$ a
Poisson group action with moment $\mu_G: M\rightarrow G^{*},$ and
$H\subset G$ a Poisson subgroup, the action of $H$ still admits a
moment map. But if $H\subset G$ is only admissible, then the
restricted action does not usually have a moment map. This is
precisely the case which is encountered in the study of the
$q$--deformed Drinfeld-Sokolov reduction.\qed
\end{rem}

\section{The $q$--deformed Drinfeld-Sokolov reduction in the case of
$SL_2$}

In this section we apply the general results of the previous section
to formulate a $q$--analogue of the Drinfeld-Sokolov reduction when
$G=SL_2$.

\subsection{Choice of $r$--matrix} Let $\g = L\sw_2$. We would like
to define a factorizable Lie bialgebra structure on $\g$ in such a way
that the resulting Poisson-Lie structure $\eta$ on $LSL_2$ and the
Poisson structure $\eta_*^q$ on $M_{2,q}$ satisfy the conditions
(ii)--(iii) of \secref{over}.

Let $\{ E,H,F \}$ be the standard basis in $\sw_2$ and $\{ E_n,H_n,F_n
\}$ be the corresponding (topological) basis of $L\sw_2 = \sw_2 \ot
\C((s))$ (here for each $A \in \sw_2$ we set $A_n = A \ot s^n \in
L\sw_2$). Let $\tau$ be the automorphism of $L\sw_2$ defined by the
formula $\tau(A(s)) = A(sq)$ (we assume that $q$ is generic). We have:
$\tau \cdot A_n = q^n A_n$. To be able to use \thmref{act1}, the
$r$--matrix $r \in L\sw_2^{\ot 2}$ defining the Lie bialgebra
structure on $L\sw_2$ has to satisfy the condition $(\tau \ot \tau)(r)
= r$. Hence the invariant bilinear form on $L\sw_2$ defined by the
symmetric part of $r$ should also be $\tau$--invariant.

The Lie algebra $L\sw_2$ has a unique (up to a non-zero constant
multiple) invariant non-degenerate bilinear form, which is invariant
under $\tau$. It is defined by the formulas
$$(E_n,F_m) = \delta_{n,-m}, \quad \quad (H_n,H_m) = 2
\delta_{n,-m},$$ with all other pairings between the basis elements
are $0$. This fixes the symmetric part of the element $r$. Another
condition on $r$ is that the subgroup $LN$ is admissible. According to
\propref{admiss}, this means that $L\n_+^\perp$ should be a Lie
subalgebra of $L\sw_2^*$.

A natural example of $r$ satisfying these two conditions is given by
the formula: \beq \label{new} r_0 = \sum_{n \in \Z} E_n \ot F_{-n} +
\frac{1}{4} H_0 \ot H_0 + \pol \sum_{n>0} H_n \ot H_{-n}.
\end{equation}
It is easy to verify that this element defines a factorizable Lie
bialgebra structure on $\g$. We remark that this Lie bialgebra
structure gives rise to Drinfeld's ``new'' realization of the
quantized enveloping algebra associated to $L\sw_2$
\cite{nr,kh-t,ked}. As we will see in the next subsection, $r_0$ can
not be used for the $q$--deformed Drinfeld-Sokolov reduction. However,
the following crucial fact will enable us to perform the
reduction. Let $L{\frak h}$ be the loop algebra of the Cartan
subalgebra ${\frak h}$ of $\sw_2$.

\begin{lem}    \label{defi}
For any $\rho \in \wedge^2 L{\frak h}$, $r_0 + \rho$ defines a
factorizable Lie bialgebra structure on $L\sw_2$, such that
$L\n_+^\perp$ is a Lie subalgebra of $L\sw_2^*$.
\end{lem}

The fact that $r_0 +\rho$ still satisfies the classical Yang-Baxter
equation is a general property of factorizable $r$--matrices
discovered in \cite{BD}. \lemref{defi} allows us to consider the
class of elements $r$ given by the formula
\beq \label{class1} r = \sum_{n \in
\Z} E_n \ot F_{-n} + \frac{1}{2} \sum_{m,n \in \Z} \phi_{n,m} \cdot
H_n \ot H_m,
\end{equation}
where $\phi_{n,m} + \phi_{m,n} = \delta_{n,-m}$.  The condition $(\tau
\ot \tau)(r) = r$ imposes the restriction $\phi_{n,m} = \phi_n
\delta_{n,-m}$, so that \eqref{class1} takes the form \beq
\label{class} r = \sum_{n \in \Z} E_n \ot F_{-n} + \frac{1}{2} \sum_{n
\in \Z} \phi_n \cdot H_n \ot H_{-n},
\end{equation}
where $\phi_n + \phi_{-n} = 1$.

\subsection{The reduction}
Recall that $M_{2,q} = LSL_2 = SL_2((s))$ consists of $2 \times 2$
matrices
\beq    \label{two}
M(s) = \begin{pmatrix} a(s) & b(s) \\ c(s) & d(s) \end{pmatrix}, \quad
\quad ad - bc = 1.
\end{equation}
We want to impose the constraint $c(s) = -1$, i.e. consider the
submanifold $M_{2,q}^J$ and take its quotient by the (free) action of
the group $$LN = \left\{ \begin{pmatrix} 1 & x(s) \\ 0 & 1
\end{pmatrix} \right\}.$$ Let $\eta$ be the Poisson-Lie structure on
$LSL_2$ induced by $r$ given by formula \eqref{class}.

Let $\eta_*^q$ be the Poisson structure on $M_{2,q}$ defined by
formula \eqref{mainb1}, corresponding to the automorphism $\tau: g(s)
\arr g(sq)$. The following is an immediate corollary of \thmref{act1},
\propref{admiss} and \lemref{defi}.

\begin{prop}
(1) The $q$--gauge action of ($LSL_2,\eta$) on ($M_{2,q},\eta^q_*$)
given by formula $g(s) \cdot M(s) = g(sq) M(s) g(s)^{-1}$ is Poisson;

(2) The subgroup $LN \subset LSL_2$ is admissible.
\end{prop}

Thus, we have satisfied properties (i) and (ii) of \secref{over}. Now
we have to choose the remaining free parameters $\phi_n$ so as to
satisfy property (iii).

The Fourier coefficients of the matrix elements of the matrix $M(s)$
given by \eqref{two} define functions on $M_{2,q}$. We will use the
notation $a_m$ for the $m$th Fourier coefficient of $a(s)$. Let
$R_{2,q}$ be the completion of the ring of polynomials in $a_m, b_m,
c_m, d_m, m \in \Z$, defined in the same way as the ring $\RN$ of
Sect.~2.3. Let $S_{2,q} \subset R_{2,q}$ be the subalgebra of
$LN$--invariant functions. Denote by $I$ be the ideal of $S_{2,q}$
generated by $\{ c_n + \delta_{n,0}, n\in\Z \}$ (the defining ideal of
$M_{2,q}^J$).

Property (iii) means that $I$ is a Poisson ideal of $s_{2,q}$, which
is equivalent to the condition that $\{ c_n,c_m \} \in I$, i.e. that
if $\{ c_n,c_m \}$ vanishes on $M_{2,q}^J$. This condition means that
the Poisson bracket of the constraint functions vanishes on the
constraint surface, i.e. the constraints are of first class according
to Dirac.

Let us compute the Poisson bracket between $c_n$'s. First, we list the
left and right gradients for the functions $a_n,b_n,c_n,d_n$ (for this
computation we only need the gradients of $c_n$'s, but we will soon
need other gradients as well). It will be convenient for us to
identify $L\sw_2$ with its dual using the bilinear form introduced in
the previous section. Note that with respect to this bilinear form the
dual basis elements to $E_n, H_n$, and $F_n$ are $F_{-n}, H_{-n}/2$,
and $E_{-n}$, respectively.

Explicit computation gives (for shorthand, we write $a$ for $a(s)$,
etc.):
$$\nb a_m = s^{-m} \begin{pmatrix} \pol a & 0 \\ c & - \pol a
\end{pmatrix}, \quad \quad \nb b_m = s^{-m} \begin{pmatrix} \pol b & 0
\\ d & - \pol b \end{pmatrix},$$
$$\nb c_m = s^{-m} \begin{pmatrix} - \pol c & a \\ 0 & \pol c
\end{pmatrix}, \quad \quad \nb d_m = s^{-m} \begin{pmatrix} - \pol d & b
\\ 0 & \pol b \end{pmatrix},$$
$$\nb' a_m = s^{-m} \begin{pmatrix} \pol a & b \\ 0 & - \pol a
\end{pmatrix}, \quad \quad \nb' b_m = s^{-m} \begin{pmatrix} - \pol b
& 0 \\ a & \pol b \end{pmatrix},$$
$$\nb' c_m = s^{-m} \begin{pmatrix} \pol c & d \\ 0 & - \pol c
\end{pmatrix}, \quad \quad \nb' d_m = s^{-m} \begin{pmatrix} - \pol d & 0
\\ c & - \pol d \end{pmatrix}.$$

Now we can compute the Poisson bracket between $c_n$'s using formula
\eqref{mainb1}: \beq \label{cbracket} \{ c_m,c_k \} = \pol \sum_{n \in
\Z} \left( \phi_n - \phi_{-n} + \phi_n q^n - \phi_{-n} q^{-n} \right)
c_{-n+m} c_{n+k}.
\end{equation}

Restricting to $M_{2,q}^J$, i.e. setting $c_n = -\delta_{n,0}$, we
obtain:
$$\{ c_m,c_k \}|_{M_{2,q}^J} = \pol \sum_{n \in \Z} \left( \phi_m -
\phi_{-m} + \phi_m q^m - \phi_{-m} q^{-m} \right) \delta_{m,-k}.$$ This
gives us the following equation on $\phi_m$'s: $$\phi_m - \phi_{-m} +
\phi_m q^m - \phi_{-m} q^{-m} = 0.$$ Together with the previous condition
$\phi_m + \phi_{-m} = 1$, this determines $\phi_m$'s uniquely:

\begin{thm}    \label{third}
{\em The Poisson structure $\eta_*^q$ satisfies property (iii) of the
$q$--deformed Drin\-feld-Sokolov reduction if and only if} $$\phi_n =
\frac{1}{1+q^n}.$$
\end{thm}

Consider the $r$--matrix \eqref{class1} with $\phi_n = (1+q^n)^{-1}$.
For this $r$--matrix, the Lie algebras defined in section
\secref{bialg} are as follows: $\g_\pm = L{\frak b}_\mp, {\frak n}_\pm
= L{\frak n}_\mp$, where ${\frak n}_+ = \C E, {\frak n}_- = \C F,
{\frak b}_\pm = {\frak h} \oplus {\frak n}_\pm$. We have:
$\g_\pm/{\frak n}_\pm \simeq L{\frak h}$. The transformation $\theta$
on $L{\frak h}$ induced by this $r$--matrix is equal to $-\tau$.

Explicitly, on the tensor product of the two $2$--dimensional
representations of $\sw_2((t))$, the $r$--matrix looks as follows:
\begin{equation}    \label{explicit}
\begin{pmatrix} \vff & 0 & 0 & 0 \\ 0 & - \vff & \delta \left(
\frac{t}{s} \right) & 0 \\ 0 & 0 & - \vff & 0 \\ 0 & 0 & 0 & \vff
\end{pmatrix},
\end{equation}
where $$\phi(x) = \pol \sum_{n \in \Z} \frac{1}{1+q^n} x^n.$$ Note
that $2 \pi \phi(xq^{1/2})$ coincides with the power series expansion
of the Jacobi elliptic function $dn$ (delta of amplitude).

Now we have satisfied all the necessary properties on the Poisson
structures and hence can perform the $q$--Drinfeld-Sokolov reduction
of \secref{qgauge} at the level of Poisson algebras. In the next
subsection we check that it indeed gives us the Poisson bracket
\eqref{vira} on the reduced space ${\cal M}_{2,q} = M_{2,q}^J/LN$.

\subsection{Explicit computation of the Poisson brackets}

Introduce the generating series $$A(z) = \sum_{n \in \Z} a_n z^{-n},$$
and the same for other matrix elements of $M(s)$ given by formula
\eqref{two}. We fix the element $r$ by setting $\phi_n = (1+q^n)^{-1}$
in formula \eqref{class} in accordance with \thmref{third}. Denote
\beq
\label{fi} \vf(z) = \sum_{n \in \Z} (\phi_n - \phi_{-n}) z^n = \sum_{n \in
\Z} \frac{1-q^n}{1+q^n} z^n.
\end{equation}
Using the formulas for the gradients of the matrix elements given in
the previous section and formula \eqref{mainb1} for the Poisson
bracket, we obtain the following explicit formulas for the Poisson
brackets.

\begin{align*}
\{ A(z),A(w) \} &= \vf \left( \frac{w}{z} \right) A(z) A(w), \\
\{ A(z),B(w) \} &= - \delta \left( \frac{w}{z} \right) A(z) B(w), \\
\{ A(z),C(w) \} &= \delta \left( \frac{w}{z} \right) A(z) C(w), \\
\{ A(z),D(w) \} &= - \vf \left( \frac{w}{z} \right) A(z) D(w),\\
\{ B(z),B(w) \} &= 0, \\
\{ B(z),C(w) \} &= \delta \left( \frac{w}{z} \right) A(z) D(w) -
\delta \left( \frac{wq}{z} \right) A(z) A(w), \\
\{ B(z),D(w) \} &= - \delta \left( \frac{wq}{z} \right) A(z) B(w), \\
\{ C(z),C(w) \} &= 0, \\
\{ C(z),D(w) \} &= \delta \left( \frac{w}{zq} \right) A(z) C(w),\\
\{ D(z),D(w) \} &= \vf \left( \frac{w}{z} \right) D(z) D(w) - \delta
\left( \frac{wq}{z} \right) C(z) B(w) + \delta \left( \frac{w}{zq}
\right) B(z) C(w).
\end{align*}

\begin{rem}
The relations above can be presented in matrix form as follows. Let
$$L(z) = \begin{pmatrix} A(z) & B(z) \\ C(z) & D(z) \end{pmatrix},$$
and consider the operators $L_1 = L \ot \on{id}, L_2 = \on{id} \ot L$
acting on $\C^2 \ot \C^2$. The $r$--matrix \eqref{explicit} also acts
on $\C^2 \ot \C^2$. Formula \eqref{mainb1} can be written as follows:
\begin{align*}
\{ L_1(z),L_2(w) \} &= \pol r^- \left( \frac{w}{z} \right) L_1(z)
L_2(w) + \pol L_1(z) L_2(w) r^- \left( \frac{w}{z} \right) \\ &-
L_1(z) r \left( \frac{wq}{z} \right) L_2(w) + L_2(w) \sigma(r) \left(
\frac{zq}{w} \right) L_1(z),
\end{align*}
where $$r^- \left( \frac{w}{z} \right) = r \left( \frac{w}{z} \right)
- \sigma(r) \left( \frac{z}{w} \right) = \begin{pmatrix} \pol \vpp & 0
& 0 & 0 \\ 0 & - \pol \vpp & \delta \left( \frac{w}{z} \right) & 0 \\
0 & - \delta \left( \frac{w}{z} \right) & - \pol \vpp & 0 \\ 0 & 0 & 0
& \pol \vpp
\end{pmatrix}.$$\qed
\end{rem}

\subsection{Reduced Poisson structure}
We know that ${\cal M}_{2,q} = M_{2,q}^J/LN$ is isomorphic to
$$\left\{
\begin{pmatrix} t(s) & 1 \\ -1 & 0 \end{pmatrix} \right\}$$ (see
Sect.~3). The ring ${\cal R}_{2,q}$ of functionals on ${\cal M}_{2,q}$
is generated by the Fourier coefficients of $t(s)$. In order to
compute the reduced Poisson bracket between them, we have to extend
them to $LN$--invariant functions on the whole $M_{2,q}$. Set \beq
\wt{t}(s) = a(s) c(sq) + d(sq) c(s).
\end{equation}
It is easy to check that the Fourier coefficients $\wt{t}_m$ of
$\wt{t}(s)$ are $LN$--invariant, and their restrictions to $M_{2,q}^J$
coincide with the corresponding Fourier coefficients of $t(s)$.

Let us compute the Poisson bracket between $\wt{t}_m$'s. Set
$$\wt{T}(z) = \sum_{m \in \Z} \wt{t}_m z^{-m}.$$ Using the explicit
formulas above, we find
\begin{align} \notag
\{ \wt{T}(z),\wt{T}(w) \} &= \vf \left( \frac{w}{z} \right) \wt{T}(z)
\wt{T}(w) \\ \label{tildet} &+ \delta \left( \frac{wq}{z} \right)
\Delta(z) c(w) c(wq^2) - \delta \left( \frac{w}{zq} \right) \Delta(w)
c(z) c(zq^2),
\end{align}
where $\Delta(z) = A(z)D(z) - B(z)C(z) = 1$. Hence, restricting to
$M_{2,q}^J$ (i.e. setting $c(z)=1$ in formula \eqref{tildet}), we
obtain:
$$\{ T(z),T(w) \} = \vf \left( \frac{w}{z} \right) T(z) T(w) + \delta
\left( \frac{wq}{z} \right) - \delta \left( \frac{w}{zq} \right).$$
This indeed coincides with formula \eqref{vira}.

\begin{rem}
Consider the subring $\wt{S}_{2,q}$ of the ring $R_{2,q}$, generated by
$c_m, \wt{t}_m, m \in \Z$. The ring $\wt{S}_{2,q}$ consists of
$LN$--invariant functionals on $M_{2,q}$ and hence it can serve as a
substitute for the ring of functions on $M_{2,q}/LN$. Let us compute
the Poisson brackets in $\wt{S}_{2,q}$. The Poisson brackets of
$\wt{t}_m$'s are given by formula \eqref{tildet}, and by construction
$\{ c_m,c_k \} = 0$. It is also easy to find that $\{ c_m,\wt{t}_k \}
= 0$. Hence $\wt{S}_{2,q}$ is a Poisson subalgebra of $R_{2,q}$. Thus, the
$q$--deformed Drinfeld-Sokolov reduction can be interpreted as
follows. The initial Poisson algebra is $R_{2,q}$. We consider its
Poisson subalgebra $\wt{S}_{2,q}$ generated by $c_m$'s and
$\wt{t}_m$'s. The ideal $I$ of $\wt{S}_{2,q}$ generated by $\{ c_m +
\delta_{m,0}, m \in \Z \}$ is a Poisson ideal. The quotient
$\wt{S}_{2,q}/I$ is isomorphic to the $q$--Virasoro algebra ${\cal
R}_{2,q}$ defined in Sect.~3.\qed
\end{rem}

\subsection{$q$--deformation of Miura transformation}

As was explained in Sect.~3.2, the $q$--Miura transformation of
\cite{FR} is the map between two (local) cross-sections of the
projection $\pi_q: M_{n,q}^J \arr M_{n,q}^J/LN$. In the case of
$LSL_2$, the first cross-section $$\left\{ \begin{pmatrix} \la(s) & 0
\\ -1 & \la(s)^{-1} \end{pmatrix} \right \}$$ is defined by the
subsidiary constraint $b(s)=0$, and the second $$\left\{
\begin{pmatrix} t(s) & 1 \\ -1 & 0 \end{pmatrix} \right
\}$$ is defined by the subsidiary constraint $d(s)=0$. The map between
them is given by the formula $$\muu_q: \la(s) \mapsto t(s) = \la(s) +
\la(sq)^{-1}.$$ Now we want to recover formula \eqref{virpois} for the
Poisson brackets between the Fourier coefficients $\la_n$ of $\la(s)$,
which makes the map $\muu_q$ Poisson.

We have already computed the Poisson bracket on the second (canonical)
cross-section from the point of view of Poisson reduction. Now we need
to compute the Poisson bracket between the functions $a_n$'s on the
first cross-section, with respect to which the map $\muu_q$ is
Poisson. This computation is essentially similar to the one outlined
in Sect.~3.2. The Poisson structure on the local cross-section is
given by the Dirac bracket, which is determined by the choice of the
subsidiary conditions, which fix the cross-section.

The Dirac bracket has the following property (see
\cite{Flato}). Suppose we are given constraints $\xi_n, n \in I$, and
subsidiary conditions $\eta_n, n \in I$, on a Poisson manifold $M$,
such that $\{ \xi_k,\xi_l \} = \{ \eta_k,\eta_l \} = 0, \forall k,l
\in I$. Let $f, g$ be two functions on $M$, such that $\{ f,\xi_k \}$
and $\{ g,\xi_k \}$ vanish on the common level surface of all $\xi_k,
\eta_k$. Then the Dirac bracket of $f$ and $g$ coincides with their
ordinary Poisson bracket.

In our case, the constraint functions are $c_m+\delta_{m,0}, m \in
\Z$, and the subsidiary conditions $b_m, m \in \Z$, which fix the
local model of the reduced space. We have: $\{ b_m,b_k \} = 0$, $\{
c_m,c_k \} = 0$, and $\{ a_m,b_k \} = 0$, if we set $b_m=0, \forall m
\in \Z$. Therefore we are in the situation described above, and the
Dirac bracket between $a_m$ and $a_k$ coincides with their ordinary
bracket. In terms of the generating function $A(z) = \sum_{m \in \Z}
a_m z^{-m}$ it is given by the formula $$\{ A(z),A(w) \} = \vf \left(
\frac{w}{z} \right) A(z) A(w),$$ which coincides with formula
\eqref{virpois}. Thus, we have proved the Poisson property of the
$q$--deformation of the Miura transformation from the point of view of
the deformed Drinfeld-Sokolov reduction.

\section{Lattice Virasoro algebra}

In this section we consider the lattice counterpart of the
Drinfeld-Sokolov reduction. Our group is thus $\GG = (SL_2)^{\Z/N\Z}$,
where $N$ is an integer, and $\tau$ is the automorphism of $G$, which
maps $(g_i)$ to $(g_{i+1})$. Poisson structures on $\GG$ which are
covariant with respect to lattice gauge transformations $x_n \mapsto
g_{n+1} x_n g_n^{-1}$ have been studied already in \cite{RIMS}
(cf. also \cite{AFS}). In order to make the reduction over the
nilpotent subgroup $\NN \subset \GG$ feasible, we have to be careful
in our choice of the $r$--matrix.

\subsection{Discrete Drinfeld-Sokolov reduction}
By analogy with the continuous case, we choose the element $r$ defining the
Lie bialgebra structure on $\g = \sw_2^{\oplus \Z/N\Z}$ as follows: $$r =
\sum_{n \in \Z/N\Z} E_n \ot F_n + \frac{1}{4} \sum_{m,n \in \Z/N\Z}
\phi_{n,m} H_n \ot H_m,$$ where $\phi_{n,m} + \phi_{m,n} = 2
\delta_{m,n}$. It is easy to see that $r$ defines a factorizable Lie
bialgebra structure on $\g$. For \thmref{act1} to be applicable, $r$ has to
satisfy the condition $(\tau \ot \tau)(r) = r$, which implies that
$\phi_{n,m} = \phi_{n-m}$.

An element of $\GG$ is an $N$--tuple $(g_i)$ of elements of $SL_2$:
$$g_k = \begin{pmatrix} a_k & b_k \\ c_k & d_k \end{pmatrix}.$$ We
consider $a_k,b_k,c_k,d_k, k \in \Z/N\Z$, as the generators of the
ring of functions on $\GG$.

The discrete analogue of the Drinfeld-Sokolov reduction consists of
taking the quotient $\MM=\GG^J/\NN$, where $\GG^J = (G^J)^{\ZN}$,
$$G^J = \left\{ \begin{pmatrix} a & b \\ -1 & d \end{pmatrix}
\right\},$$ and $\NN = N^{\ZN}$, acting on $G^J$ by the formula \beq
\label{dgauge} (h_i) \cdot (g_i) = (h_{i+1} g_i h_i^{-1}).
\end{equation}
It is easy to see that $$\MM \simeq \left \{ \begin{pmatrix} t_i &
1 \\ -1 & 0 \end{pmatrix}_{i \in \ZN} \right\}.$$

The element $r$ with $\phi_{n,m} = \phi_{n-m}, \phi_k + \phi_{-k} =
2\delta_{k,0}$, defines a Lie bialgebra structure on $\g$ and Poisson
structures $\eta, \eta_*^\tau$ on $\GG$. According to \thmref{act1},
the action of ($\GG,\eta$) on ($\GG,\eta_*^\tau$) given by formula
\eqref{dgauge} is Poisson.

As in the continuous case, for the Poisson structure $\eta_*^\tau$ to
be compatible with the discrete Drinfeld-Sokolov reduction, we must
have:
\beq    \label{vanish}
\{ c_n,c_m \}|_{\GG^J} = 0.
\end{equation}
Explicit calculation analogous to the one made in the previous
subsection shows that \eqref{vanish} holds if and only if $$\phi_{n-1}
+ 2\phi_n + \phi_{n+1} = 2 \delta_{n,0} + 2 \delta_{n+1,0}.$$ The
initial condition $\phi_0 = 1$ and periodicity condition give us a
unique solution: for odd $N$, $\phi_k = (-1)^k$; for even $N$, $\phi_k
= (-1)^k \left( 1 - \dfrac{2k}{N} \right)$. In what follows we
restrict ourselves to the case of odd $N$ (note that in this case the
linear operator $\on{id}+\tau$ is invertible).

Continuing as in the previous subsection, we define $$\wt{t}_n = a_n
c_{n+1} + d_{n+1} c_n, \quad \quad n \in \ZN.$$ These are
$\NN$--invariant functions on $\GG$. We find in the same way as in the
continuous case:
\beq    \label{wtt}
\{ \wt{t}_n,\wt{t}_m \} = \vf_{n-m} \wt{t}_n \wt{t}_m + \delta_{n,m+1}
c_m c_{m+2} - \delta_{n+1,m} c_n c_{n+2},
\end{equation}
$$\{ \wt{t}_n,c_m \} = 0, \quad \quad \{ c_n,c_m \} = 0,$$ where
$$\vf_k = \pol (\phi_k - \phi_{-k}) = \left\{ \begin{array}{cc} 0, &
k=0, \\ (-1)^k, & k\neq 0.
\end{array} \right.$$
The discrete Virasoro algebra $\C[t_i]_{i \in \ZN}$ is the quotient of
the Poisson algebra $\C[\wt{t}_i,c_i]_{i \in \ZN}$ by its Poisson
ideal generated by $c_i+1, i \in \ZN$. From formula \eqref{wtt} we
obtain the following Poisson bracket between the generators $t_i$:
\beq
\label{dvir} \{ t_n,t_m \} = \vf_{n-m} t_n t_m + \delta_{n,m+1} -
\delta_{n+1,m}.
\end{equation}

The discrete Miura transformation is the map from the local
cross-section $$\left\{ \begin{pmatrix} \la_n & 0 \\ -1 & \la_n^{-1}
\end{pmatrix} \right \}$$ to $\MM$,
\beq    \label{dmiura}
\la_n \mapsto t_n = \la_n + \la_{n+1}^{-1}.
\end{equation}
It defines a Poisson map $\C[\la_i^\pm]_{i \in \ZN} \arr \C[t_i]_{i
\in \ZN}$, where the Poisson structure on the latter is given by the
formula \beq \label{dheis} \{ \la_n,\la_m \} = \vf_{n-m} \la_n \la_m.
\end{equation}

\begin{rem}
The Poisson algebra $\C[t_i]_{i \in \ZN}$ can be considered as a
regularized version of the $q$--deformed Virasoro algebra when
$q=\ep$, where $\ep$ is a primitive $N$th root of unity. Indeed, we
can then consider $t(\ep^i), i \in \ZN$, as generators and truncate in
all power series appearing in the relations, summations over $\Z$ to
summations over $\ZN$ divided by $N$. This means that we replace
$\phi(\ep^n)$ given by formula \eqref{vira} by
$$\wt{\phi}(\ep^n) = \frac{1}{N} \sum_{i \in \ZN}
\frac{1-\ep^i}{1+\ep^i} \ep^{ni},$$ and $\delta(\ep^n)$ by
$\delta_{n,0}$. The formula for the Poisson bracket then becomes: $$\{
t(\ep^n),t(\ep^m) \} = \wt{\phi}(\ep^{m-n}) t(\ep^n) t(\ep^m) +
\delta_{n,m+1} - \delta_{n+1,m}.$$ If we set $t(\ep^i)=t_i$, we
recover the Poisson bracket \eqref{dvir}, since it is easy to check
that $\wt{\phi}(\ep^{m-n}) = \vf_{n-m}$.

One can apply the same procedure to the $q$--deformed $\W$--algebras
associated to $\sw_N$ and obtain lattice Poisson algebras. It would be
interesting to see whether they are related to the lattice
$\W$--algebras studied in the literature, e.g., in \cite{Be,Bo}. In
the case of $\sw_2$, this connection is described in the next
subsection.\qed
\end{rem}

\subsection{Connection with Faddeev-Takhtajan-Volkov algebra}
The Poisson structures \eqref{dvir} and \eqref{dheis} are nonlocal,
i.e. the Poisson brackets between distant neighbors on the lattice
are nonzero. However, one can define closely connected Poisson
algebras possessing local Poisson brackets; these Poisson algebras can
actually be identified with those studied by L.~Faddeev,
L.~Takhtajan, and A.~Volkov.

Let us first recall some results of \cite{FR} concerning the
continuous case. As was explained in \cite{FR}, one can associate a
generating series of elements of the $q$--Virasoro algebra to an
arbitrary finite-dimensional representation of $\sw_2$. The series
$T(z)$ considered in this paper corresponds to the two-dimensional
representation. Let $T^{(2)}(z)$ be the series corresponding to the
three-dimensional irreducible representation of $sw_2$. We have the
following identity \cite{FR}
$$T(z) T(zq) = T^{(2)}(z) + 1,$$ which can be taken as the definition of
$T^{(2)}(z)$. From formula \eqref{virmiura} we obtain:
\begin{align*}
T^{(2)}(z) &= \La(z) \La(zq) + \La(z) \La(zq^2)^{-1} + \La(zq)^{-1}
\La(zq^2)^{-1} \\ &= A(z) + A(z) A(zq)^{-1} + A(zq)^{-1},
\end{align*}
where
\beq    \label{az}
A(z) = \La(z) \La(zq)
\end{equation}
(the series $A(z)$ was introduced in Sect.~7 of \cite{FR}). From
formula \eqref{virpois} we find:
$$\{ A(z),A(w) \} = \left( \delta \left( \frac{w}{zq} \right) - \delta
\left( \frac{wq}{z} \right) \right) A(z) A(w).$$ It is also easy to find
\begin{align*}    \label{three}
\{ T^{(2)}(z),T^{(2)}(w) \} &= \left( \delta \left( \frac{w}{zq}
\right) - \delta \left( \frac{wq}{z} \right) \right) \left( T^{(2)}(z)
T^{(2)}(w) - 1 \right) \\ &+ \delta \left( \frac{wq^2}{z} \right) T(w)
T(wq^3) - \delta \left( \frac{w}{zq^2} \right) T(z) T(zq^3).
\end{align*}

We can use the same idea in the lattice case. Let $\nu_n = \la_n
\la_{n+1}$; this is the analogue of $A(z)$. We have: \beq \label{nun}
\{ \nu_n,\nu_m \} = (\delta_{n+1,m} - \delta_{n,m+1}) \nu_n \nu_m,
\end{equation}
and hence $\C[\nu_i^\pm]$ is a Poisson subalgebra of $\C[\la_n^\pm]$
with local Poisson brackets. We can also define $t^{(2)}_n = t_n
t_{n+1} - 1$. The Poisson bracket of $t^{(2)}_n$'s is local:
\begin{align}    \label{t2}
\{ t^{(2)}_n,t^{(2)}_m \} &= \left( \delta_{n+1,m} - \delta_{n,m+1}
\right) \left( t^{(2)}_n t^{(2)}_m - 1 \right) \\ \notag &+
\delta_{n,m+2} t_m t_{m+3} - \delta_{n+2,m} t_n t_{n+3}.
\end{align}
Unfortunately, it does not close on $t^{(2)}_n$'s, so that
$\C[t^{(2)}_i]$ is not a Poisson subalgebra of $\C[t_i]$. But let us
define formally
\beq    \label{sn}
s_n = \frac{1}{1+t^{(2)}_n} = t_n^{-1} t_{n+1}^{-1} =
\frac{1}{(1+\nu_n)(1+\nu_{n+1}^{-1})}.
\end{equation}
Then from formulas \eqref{sn} and \eqref{nun} we find:
\begin{align}    \label{fad}
\{ s_n,s_m \} = & s_n s_m \big( (\delta_{n+1,m} - \delta_{n,m+1})(1 -
s_n - s_m) - \\ \notag & - s_{n+1} \delta_{n+2,m} + s_{m+1}
\delta_{n,m+2} \big).
\end{align}
Thus, the Poisson bracket closes among $s_n$'s and defines a Poisson
structure on $\C[s_i]_{i\in\ZN}$.

The Poisson algebra $\C[s_i]_{i\in\ZN}$ with Poisson bracket
\eqref{fad} was first introduced by Faddeev and Takhtajan in \cite{ft}
(see formula (54)). We see that it is connected with our version of
the discrete Virasoro algebra, $\C[t_i]$, by a change of variables
\eqref{sn}. The Poisson algebra $\C[\nu_i^\pm]$ and the Poisson map
$\C[\nu_i^\pm] \arr \C[s_n]$ given by formula \eqref{sn} were
introduced by Volkov in \cite{v1} (see formulas (2) and (23))
following \cite{ft}; see also related papers \cite{v2,fv}. This map is
connected with our version \eqref{dmiura} of the discrete Miura
transformation by a change of variables.


\begin{thebibliography}{99}

\bibitem{A} Adler, M., On a trace functional for formal
pseudodifferential operators and the symplectic structure of the
Korteweg--de Vries type equations, {\em Inv. Math.} {\bf 50} (1979),
219-248.

\bibitem{AFS} Alekseev, A., Faddeev, L., Semenov-Tian-Shansky, M.,
Hidden quantum group inside Kac--Moody algebras, {\em
Commun. Math. Phys.} {\bf 149} (1992), 335-345.

\bibitem{Miwa} Asai, Y., Jimbo, M., Miwa, T., Pugay, Ya., Bosonization of
vertex operators for the $A^{(1)}_{n-1}$ vertex models, J. Phys. {\bf
A29} (1996), 6595-6616.

\bibitem{AKOS} Awata, H., Kubo, H., Odake, S., Shiraishi, J., Quantum
$W_n$ algebras and Macdonald polynomials, {\em Comm. Math. Phys.} {\bf
179} (1996), 401-416.

\bibitem{BD} Belavin, A.A., Drinfeld V.G., Solutions of the classical
Yang-Baxter equation for simple Lie algebras, {\em
Funct. Anal. Appl.,} {\bf 16} (1981), 159-80.

\bibitem{Be} Belov, A.D., Chaltikian, K.D., Lattice analogues of
$W$--algebras and classical integrable equations, {\em Phys. Lett.}
(1993), 268-274.

\bibitem{Bo} Bonora, L., Colatto, L.P., Constantinidis, C.P., Toda
lattice field theories, discrete W--algebras, Toda lattice hierarchies
and quantum groups, {\em Phys. Lett.} {\bf B387} (1996), 759.

\bibitem{Hou} Ding, X.-M., Hou, B.-Y., Zhao, L., $\hbar$-(Yangian)
deformation of Miura map and Virasoro algebra, Preprint q-alg/9701014.

\bibitem{Dr} Drinfeld, V.G., Hamiltonian structures on Lie groups,
Lie bialgebras and the geometric meaning of the classical Yang--Baxter
equation, {\em Sov. Math. Dokl.} {\bf 27} (1983), 68-71.

\bibitem{nr} Drinfeld, V.G., A new realization of Yangians and
quantized affine algebras, {\em Sov. Math. Dokl.} {\bf 36} (1988),
212-216.

\bibitem{DS} Drinfeld, V.G., Sokolov, V.V., Lie algebras and equations
of Korteweg-de Vries type, {\em Sov. Math. Dokl.} {\bf 23} (1981),
457-62; {\em J. Sov. Math.} {\bf 30} (1985), 1975-2035.

\bibitem{GD} Gelfand, I.M., Dickey, L.A., Family of Hamiltonian
structures connected with integrable nonlinear equations. Collected
papers of I.M.  Gelfand, vol. 1, Springer-Verlag (1987), 625-46.

\bibitem{fv} Faddeev, L.D., Volkov, A.Yu., Abelian current algebras
and the Virasoro algebra on the lattice, {\em Phys. Lett.} {\bf B315}
(1993), 311-8.

\bibitem{FF}  Feigin, B., Frenkel, E., Affine Lie algebras at the critical
level and Gelfand-Dikii algebras, {\em Int. J. Math. Phys.} {\bf A7}, suppl.
A1 (1992), 197-215.

\bibitem{ell} Feigin, B., Frenkel, E., Quantum ${\cal W}$-algebras and
elliptic algebras, {\em Comm. Math. Phys.} {\bf 178} (1996), 653-678;
q-alg/9508009.

\bibitem{Flato} Flato, M., Lichnerowicz, A., Sternheimer, D.,
Deformation of Poisson brackets, Dirac brackets and applications, {\em
J. Math. Phys.} {\bf 17} (1976), 1754.

\bibitem{Fr} Frenkel, E., Deformations of the KdV hierarchy and
related soliton equations, {\em Int. Math. Res. Notices} {\bf 2}
(1996) 55-76; q-alg/9511003.

\bibitem{EF} Frenkel, E., Affine Kac-Moody algebras at the critical
level and quantum Drinfeld-Sokolov reduction. PhD Thesis, Harvard
University, 1991.

\bibitem{FR} Frenkel, E., Reshetikhin, N., Quantum affine algebras and
deformations of the Virasoro algebra and ${\cal W}$-algebras, {\em
Comm. Math. Phys.} {\bf 178} (1996), 237-264; q-alg/9505025.

\bibitem{Kac} Kac, V., Infinite-dimensional Lie algebras, Third
edition, Cambridge University Press, 1990.

\bibitem{ked} Kedem, R., Singular $R$--matrices and Drinfeld's
comultiplication, Preprint q-alg/9611001.

\bibitem{kh-t} Khoroshkin, S.M., Tolstoy, V.N., On Drinfeld's
realization of quantum affine algebras, {\em J. Geom. Phys.} {\bf 11}
(1993), 445-52.

\bibitem{Ko} Kostant, B., On Whittaker vectors and representation
theory, {\em Invent. Math.} {\bf 48} (1978), 101-184.

\bibitem{Lu}  Lu, J.H., Momentum mapping and reduction of Poisson actions. 
In: Symplectic geometry, groupoids and integrable systems, Berkeley,
1989. P.Dazord and A.Weinstein (eds.), pp. 209-226. Springer-Verlag.

\bibitem{LP} Lukyanov, S., Pugay, Ya., Multi-point local height
probabilities in the integrable RSOS models, {\em Nucl. Phys.} {\bf
B473} (1996), 631.

\bibitem{ext} Reshetikhin, N.Yu., Semenov-Tian-Shansky, M.A., Central
extensions of quantum current groups, {\em Lett. Math. Phys.} {\bf 19}
(1990), 133-42.

\bibitem{fact} Reshetikhin, N.Yu., Semenov-Tian-Shansky, M.A., Quantum
R-matrices and factorization problems. In: Geometry and Physics,
essays in honor of I.M.Gelfand, S.Gindikin and I.M.Singer (eds.),
pp. 533-550. North Holland, Amsterdam--London--New York, 1991.

\bibitem{rmatr} Semenov-Tian-Shansky, M.A., What is a classical
$r$--matrix, {\em Funct. Anal. Appl.} {\bf 17} (1983), 17-33.

\bibitem{RIMS}  Semenov-Tian-Shansky M.A., Dressing action transformations
and Poisson--Lie group actions. {\em Publ. RIMS.} {\bf 21} (1985),
1237-1260.

\bibitem{dual} Semenov-Tian-Shansky, M.A., Poisson Lie groups, quantum
duality principle and the quantum double. {\em Contemporary Math.}
{\bf 175}, 219-248.

\bibitem{SS} Semenov-Tian-Shansky M.A., Sevostyanov A.V.,
Drinfeld-Sokolov reduction for difference operators and deformations
of ${\cal W}$--algebras II. General semisimple case, Preprint
q-alg/9702016.

\bibitem{SKAO} Shiraishi, J., Kubo, H., Awata, H., Odake, S., A quantum
deformation of the Virasoro algebra and the Macdonald symmetric
functions, {\em Lett. Math. Phys.} {\bf 38} (1996), 33-51.

\bibitem{St} Steinberg, R., Regular elements of semisimple algebraic Lie
groups, {\em Publ. Math. I.H.E.S.}, {\bf 25} (1965), 49-80.

\bibitem{ft} Takhtajan, L.A., Faddeev, L.D., Liouville model on the
lattice, {\em Lect. Notes in Phys.} {\bf 246} (1986), 166-179.

\bibitem{v1} Volkov, A.Yu., Miura transformation on a lattice,
{\em Theor. Math. Phys.} {\bf 74} (1988), 96-99.

\bibitem{v2} Volkov, A.Yu., Quantum Volterra model, Preprint
HU-TFT-92-6 (1992).

\bibitem{W} Weinstein, A., Local structure of Poisson manifolds, {\em
J.  Different. Geom.} {\bf 18} (1983), 523-558.
\end{thebibliography}
\end{document}